\documentclass{article}
\usepackage{fullpage}
\usepackage[utf8]{inputenc}
\usepackage{graphicx}
\usepackage{amsmath, amssymb}

\title{On the Generalization of  Rademacher's Proof to $\displaystyle \vartheta_1(z,\tau)$}
\author{Maher Me'meh \& Ali Saraeb }
\date{}

\usepackage{titlesec}
\titleformat*{\section}{\centering\large\bfseries}
\titlespacing*{\section}{0pt}{\baselineskip}{\baselineskip}
\begin{document}

\maketitle
\noindent
\textbf{Abstract}: In this paper, we generalize  Rademacher's proof of the transformation law of the eta function to the Jacobi theta function using Residue calculus. 
\section*{
    INTRODUCTION}
Jacobi theta functions, which are quasi-periodic entire functions closely related to modular and Jacobi forms, have been studied extensively, since the beginning of the $19^{th}$ century by mathematicians like Jacobi, Weirstrass and Eisenstein,  for the relations they hold with fields like quantum theory, PDEs, Moduli spaces, quadratic forms and some important topics like the integral representation of the Riemann zeta function and in general $L$-functions.
\vspace{0.25cm}
\\
We define $\vartheta_1(z,\tau)$ by
\begin{align} \vartheta_1(z,\tau)=-i\sum_{n=-\infty}^{\infty} (-1)^n q^{(n+1/2)^2}e^{(2n+1)i\pi z} 
\end{align} 
where $z\in ~\mathbf{C}$ and $q=e^{2\pi i\tau}$ where $\tau\in \mathbf{H}$, the upper half plane.\\
However, using the well-known Jacobi triple product identity with $w=e^{\pi i z}$, we have
\begin{align}
\displaystyle  \sum_{n=-\infty}^{+\infty} w^{2n}q^{n^2} = \prod_{m=1}^{\infty} (1-q^{2m})(1+ w^2 q^{2m-1})(1+ w^{-2} q^{2m-1}).
\end{align}
Having applied an appropriate change of variable to (2), we arrive at an exquisite product representation of $\vartheta_1$ as shown below
\begin{align}
\vartheta_1(z,\tau) = -i w q^{1/4} \prod_{n=1}^{\infty} (1-q^{2n})(1- w^2 q^{2n})(1- w^{-2} q^{2n-2}).
 \end{align} 
It is important to note that the Dedekind eta function, defined by
\begin{align}
    \eta(\tau)=e^{\frac{\pi i \tau}{12}}\prod_{n=1}^{\infty}(1-e^{2\pi i n\tau})
\end{align}
is related to the theta functions mentioned under their transformation law over an element \\  $A=\left(\begin{array}{cc}
  a   &  b\\
   c  & d\\
\end{array}\right)\in$ $\Gamma=$SL$_2(\mathbf{Z})=\left\lbrace \left(\begin{array}{cc}
   a  &  b\\
   c  &d \\
\end{array}\right) \left. \right| ~a,b,c,d \in \mathbf{Z},~ ad-bc=1\right\rbrace $. Such relations have to do with the multiplier system $\varepsilon(A)$ that appears in the transformation law of the eta function under a matrix $A\in$ SL$_2$(Z)
\[\eta(A\tau)=\varepsilon(A)(-i(c\tau+d))^{1/2}\eta(\tau).\]
Another connection is noted that most of the proofs given for the transformation law of the eta function can be generalized for the Jacobi theta function such as Iseki's and Gordon's proofs as found in \cite{g. iseki} and \cite{g. gordon}.
\vspace{0.2cm}
\\
Seigel, in \cite{Ap}, gave an elegant proof for the transformation law of the eta function under the inversion $S=\left(\begin{array}{cc}
  0   &  -1\\
   1  & 0\\
\end{array}\right)$ using residue calculus on the function
\[F_N=-\frac{1}{8z}\mathrm{\cot(\pi i Nz) \, \cot \left(\frac{\pi Nz}{y} \right)}.\]
Several extensions of Seigel's proof have been introduced, such as Raji's in \cite{Raji} and Rademecher's in \cite{Rd}.  Rademacher generalized Seigel's proof for the general transformation law of the eta function using the following function 
\begin{align}
    F_N=\frac{-1}{4ix} \cot(\pi Nx) \, \cot \left(\frac{\pi Nx}{z} \right)+\sum_{\mu=1}^{k-1}\frac{1}{x} \cdot \frac{e^{2\pi N\mu x/k}}{1-e^{2\pi Nx}} \cdot \frac{e^{2\pi i N\mu^* x/kz}}{1-e^{2\pi i N x/z}},
\end{align}
where the aforementioned symbols in (5) are to be understood from the context in $\cite{Rd}$.\\ \\
Jacobi forms, as demonstrated by Dabholkar, Murthy, and Zagier in their book $\cite{dobholkar2012quantum}$, have certain connections to Jacobi theta functions.
As found in section 4, a Jacobi form with index 1 and weight 10 denoted by \(\varphi_{10,1}(\tau, z)\), can be expressed by:

\[
\varphi_{10,1}(\tau, z) = \eta^{18}(\tau) \, \vartheta^2(\tau,z).
\]
Also for weight -2 and index 1, the Jacobi form is representable by
\[\varphi_{-2,1}(\tau,z)=\frac{\vartheta^2(\tau,z)}{\eta^6(\tau)}.\]
Other indirect connections between $\vartheta_1(\tau,z)$ and Dedekind eta function to Jacobi forms which use relations between $\vartheta_1(\tau,z)$ and other Jacobi theta functions like $\vartheta_2,~\vartheta_3,~\vartheta_4$, highlight the relevance and significance of providing an extension of  Rademacher's proof for $\eta$ to $\vartheta_1$, in hope that it can be generalized. \\ \\
Although some proofs have been provided for the transformation laws of the theta functions using the theory of quadratic forms, (see \cite{Rich} and  \cite{Yurii}-6.6), we intend to establish a newer and an easier proof for one of the transformation laws by appealing to residue calculus. Thus, in this paper we generalize  Rademacher's proof for the eta function by devising a special function ``$F_N$", to prove the transformation law of the Jacobi theta function $\vartheta_1(\tau,z)$.\\ \\ 
We now present the theorem of the transformation law of $\vartheta_1(z,\tau)$ under the action of the full modular group, $\Gamma$.
\vspace{0.3cm}
\\
$\textbf{Theorem }$ For $A=\left(\begin{array}{cc}
   a  &  b\\
   c  &d \\
\end{array}\right)~\in~\Gamma$ such that $c>0$, $\tau\in \mathbf{H}$, and $z\in \mathbf{C}$, the transformation law of $\vartheta_1$ is given by
\begin{equation} \displaystyle
\vartheta_1\left(\frac{z}{c\tau+d},\frac{a\tau+b}{c\tau+d}\right)=\varepsilon_1(A)\left(-i(c\tau+d)\right)^{1/2}e^{\frac{\pi i c z^2}{c\tau+d}}\vartheta_1(z,\tau),
\end{equation}
where 
 \[\varepsilon_1(A)=-i\varepsilon^3(A)\]
Here $\varepsilon$ appears in the transformation law of the Dedekind eta function as mentioned in the introduction and is defined by such:
\[\varepsilon(A)=\exp \left(\pi i \left(\frac{a+d}{12c}+s(-d,c) \right) \right),\]
where
\[\displaystyle s(h,k)=\sum_{r=1}^{k-1}\frac{r}{k}\left(\frac{hr}{k}-\left[\frac{hr}{k}\right]-\frac{1}{2}\right)\] is the Dedekind sum for $k>0$ and $(k,h)=1$.  
\section*{
    PROOF OF THE MAIN THEOREM}
If $z \in \mathbb{C} \backslash \{0\}$ and $x \in \mathbb{C},$ let $z^x = e^{x \\log z}$, where $-\pi < \arg z \leq \pi.$ \\ \\
Taking logarithm both sides, we observe that proving (6) is equivalent to proving 
\begin{align}
\log \left(\vartheta_1\left(\frac{z}{c\tau+d},\frac{a\tau+b}{c\tau+d}\right)\right)&= \log\left(\vartheta_1\left(z,\tau\right)\right) +\log(\epsilon_1(A)) +\frac{\pi i c z^2}{c\tau+d} +\frac{1}{2}\log \left(-i(c\tau+d) \right) \nonumber \\ &=  \log\left(\vartheta_1\left(z,\tau\right)\right) - \frac{\pi i}{2} + 3\pi i \left(\frac{a+d}{12c}+s(-d,c) \right)  +\frac{\pi i c z^2}{c\tau+d}  \nonumber \\  ~~~~~~~& ~~~+\frac{1}{2}\log(-i(c\tau+d)).
\end{align}
At this stage, we introduce a classical change of variable
\[-i(c\tau+d)=v~~~~~~~~a=H,~ c=k,~and~h=-d,\]
such that $ k>0, ~(h,k)=1,$ and $Hh\equiv -1 \pmod{k}.$ 
\\
\\
In virtue of this change of variable, we have
\[\tau=\frac{iv+h}{k}~~~\&~~~~\frac{a\tau+b}{c\tau+d}=\frac{1}{k}\left(H+\frac{i}{v}\right).\] 
Thus (7) becomes \\ 
\begin{equation}
\log\left(\vartheta_1\left(\frac{z}{iv},\frac{1}{k}\left(H+\frac{i}{v}\right)\right)\right)\\= - \frac{\pi i}{2} +  \frac{H-h}{4k}+ 3\pi is(h,k) + \frac{1}{2}\log(v)+\frac{\pi i k z^2}{iv}+\log\left(\vartheta_1\left(z,\frac{iv+h}{k}\right)\right).
\end{equation}
In other words, we have to show 
\begin{equation}
\log\left(\vartheta_1\left(z,\frac{iv+h}{k}\right)\right)-\log\left(\vartheta_1\left(\frac{z}{iv},\frac{1}{k}\left(H+\frac{i}{v}\right)\right)\right) - \frac{\pi i}{2} +  \frac{H-h}{4k}+ 3\pi is(h,k)+\frac{\pi i k z^2}{iv} = - \frac{1}{2}\log(v).
\end{equation}
Let's first consider the term $\displaystyle \log \left(\vartheta_1\left(z,\frac{iv+h}{k}\right)\right).$
\\ \\
First, note that 
\begin{align}
\log\left(\vartheta_1\left(z,\tau\right)\right) &=  \sum_{n=1}^{\infty} \left[\log \left(1-e^{2n\pi i \tau} \right)+ \log \left(1-e^{2\pi i z}e^{2n\pi i \tau} \right)  +\log \left(1-e^{-2\pi i z}e^{2(n-1)\pi i \tau} \right)\right] - \frac{\pi i}{2} +\pi i z+ \frac{\pi i \tau}{4} 
\end{align}
Therefore, we have 
\begin{align}
\log\left(\vartheta_1\left(z,\frac{iv+h}{k}\right)\right) &=\sum_{n=1}^{\infty} \log \left(1-e^{\frac{2n\pi i}{k}(h+iv)} \right)+ \sum_{n=1}^{\infty} \log \left(1-e^{2 \pi i z}e^{\frac{2n\pi i}{k}(h+iv)} \right) \nonumber \\ &~~~+ \sum_{n=1}^{\infty}  \log \left(1-e^{-2 \pi i z}e^{\frac{2(n-1)\pi i}{k}(h+iv)} \right)  - \frac{\pi i}{2} +\pi i z+ \frac{\pi i (iv+h) }{4k} \nonumber \\ 
&=\sum_{\mu=1}^{k}\sum_{q=0}^{\infty} \log \left(1-e^{\frac{2\pi i  h \mu}{k}} e^{\frac{-2\pi v  \left(qk + \mu \right)}{k}}\right) 
+ \sum_{\mu=1}^{k}\sum_{q=0}^{\infty} \log \left(1-e^{2 \pi i z} e^{\frac{2\pi i  h \mu}{k}} e^{\frac{-2\pi v  \left(qk + \mu \right)}{k}}\right)
\nonumber \\ &~~~+ \sum_{\mu=1}^{k}\sum_{q=0}^{\infty} \log \left(1-e^{-2 \pi i z} e^{\frac{2\pi i  h (\mu-1)}{k}} e^{\frac{-2\pi v \left(qk + \mu -1\right)}{k}}\right) - \frac{\pi i}{2} +\pi i z+ \frac{\pi i (iv+h) }{4k} \nonumber \\
&= - \sum_{\mu=1}^{k}\sum_{q=0}^{\infty}\sum_{n=1}^{\infty} \frac{1}{n} e^{\frac{2\pi i n h  \mu}{k}} e^{\frac{-2\pi n v  \left(qk + \mu \right)}{k}} 
-\sum_{\mu=1}^{k}\sum_{q=0}^{\infty}\sum_{n=1}^{\infty} \frac{1}{n} e^{2 \pi i n z} e^{\frac{2\pi i n  h \mu}{k}} e^{\frac{-2\pi n v  \left(qk + \mu \right)}{k}}
\nonumber \\ & ~~~- \sum_{\mu=1}^{k}\sum_{q=0}^{\infty}\sum_{n=1}^{\infty} \frac{1}{n} e^{-2 \pi i n z} e^{\frac{2\pi i n  h (\mu-1)}{k}} e^{\frac{-2\pi n v \left(qk + \mu -1\right)}{k}} - \frac{\pi i}{2} +\pi i z+ \frac{\pi i (iv+h) }{4k} \end{align}
Thus, we have  
\begin{align}
\log\left(\vartheta_1\left(z,\frac{iv+h}{k}\right)\right)  &= - \sum_{\mu=1}^{k}\sum_{n=1}^{\infty} \frac{1}{n} e^{\frac{2\pi i n h  \mu}{k}} \frac{e^{\frac{-2\pi n v  \mu}{k}}}{1- e^{-2 \pi n v}}
-\sum_{\mu=1}^{k}\sum_{n=1}^{\infty} \frac{1}{n} e^{2 \pi i n z} e^{\frac{2\pi i n h  \mu}{k}} \frac{e^{\frac{-2\pi n v  \mu}{k}}}{1- e^{-2 \pi n v}} 
\nonumber \\ &~~~- \sum_{\mu=1}^{k}\sum_{n=1}^{\infty} \frac{1}{n} e^{-2 \pi i n z} e^{\frac{2\pi i n  h (\mu-1)}{k}} \frac{e^{\frac{-2\pi n v \left( \mu -1\right)}{k}} }{1- e^{-2 \pi n v}}
- \frac{\pi i}{2} +\pi i z+ \frac{\pi i (iv+h) }{4k} \end{align}
Analogously to (12), we get 
\begin{align}
\log\left(\vartheta_1\left(\frac{z}{iv},\frac{1}{k}\left(H+\frac{i}{v}\right)\right)\right) 
&= - \sum_{w=1}^{k}\sum_{n=1}^{\infty} \frac{1}{n} e^{\frac{2\pi i n H  w}{k}} \frac{e^{\frac{-2\pi n   w}{kv}}}{1- e^{ \frac{-2 \pi n}{v}}}
-\sum_{w=1}^{k}\sum_{n=1}^{\infty} \frac{1}{n} e^{2 \pi  n \frac{z}{v}} e^{\frac{2\pi i n H  w}{k}} \frac{e^{\frac{-2\pi n   w}{kv}}}{1- e^{\frac{-2 \pi n}{v}}} 
\nonumber \\ &~~~- \sum_{w=1}^{k}\sum_{n=1}^{\infty} \frac{1}{n} e^{-2 \pi  n \frac{z}{v}} e^{\frac{2\pi i n  H (w-1)}{k}} \frac{e^{\frac{-2\pi n  \left( w -1\right)}{kv}} }{1- e^{\frac{-2 \pi n}{v}}} - \frac{\pi i}{2} +\pi \frac{z}{v}+ \frac{\pi i (i/v+H) }{4k} \end{align}
\begin{center}
\begin{figure}
    \centering
    \includegraphics[scale=0.34]{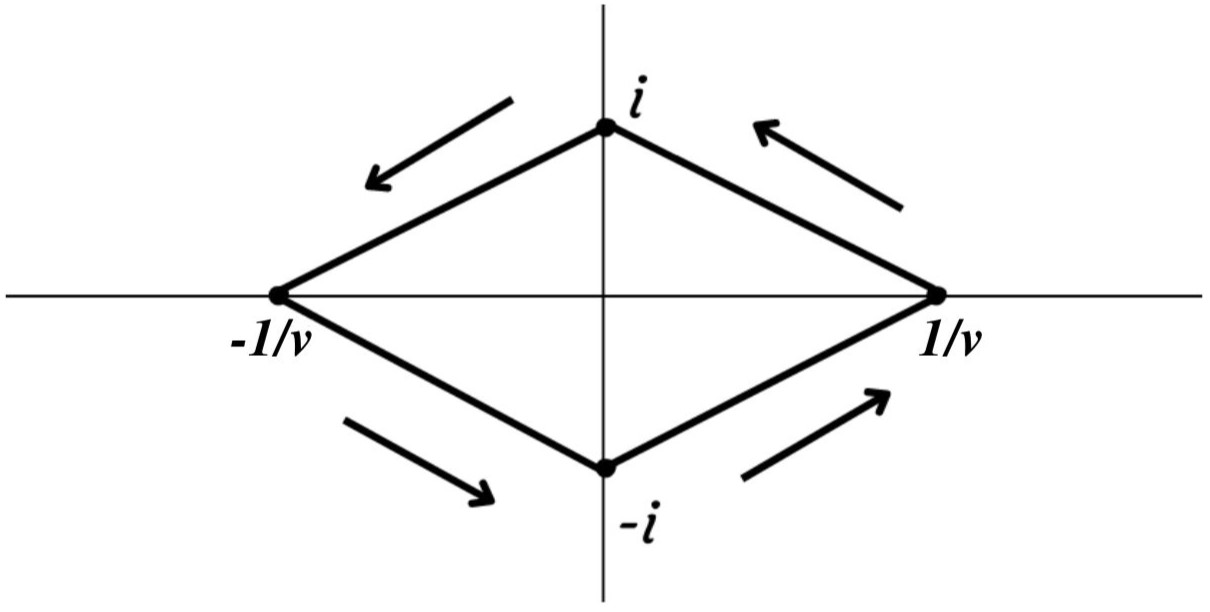}
    \caption{ Contour C}
\end{figure}
\end{center}
Plugging (12) and (13) in (9), what is required to prove becomes
\begin{align}
&\sum_{w=1}^{k}\sum_{n=1}^{\infty} \frac{1}{n} e^{\frac{2\pi i n H  w}{k}} \frac{e^{\frac{-2\pi n   w}{kv}}}{1- e^{ \frac{-2 \pi n}{v}}} +\sum_{w=1}^{k}\sum_{n=1}^{\infty} \frac{1}{n} e^{2 \pi  n \frac{z}{v}} e^{\frac{2\pi i n H  w}{k}} \frac{e^{\frac{-2\pi n   w}{kv}}}{1- e^{\frac{-2 \pi n}{v}}} 
\nonumber  \\ &+ \sum_{w=1}^{k}\sum_{n=1}^{\infty} \frac{1}{n} e^{-2 \pi  n \frac{z}{v}} e^{\frac{2\pi i n  H (w-1)}{k}} \frac{e^{\frac{-2\pi n  \left( w -1\right)}{kv}} }{1- e^{\frac{-2 \pi n}{v}}} - \sum_{\mu=1}^{k}\sum_{n=1}^{\infty} \frac{1}{n} e^{\frac{2\pi i n h  \mu}{k}} \frac{e^{\frac{-2\pi n v  \mu}{k}}}{1- e^{-2 \pi n v}}
\nonumber \\ 
 &-\sum_{\mu=1}^{k}\sum_{n=1}^{\infty} \frac{1}{n} e^{2 \pi i n z} e^{\frac{2\pi i n h  \mu}{k}} \frac{e^{\frac{-2\pi n v  \mu}{k}}}{1- e^{-2 \pi n v}} 
 - \sum_{\mu=1}^{k}\sum_{n=1}^{\infty} \frac{1}{n} e^{-2 \pi i n z} e^{\frac{2\pi i n  h (\mu-1)}{k}} \frac{e^{\frac{-2\pi n v \left( \mu -1\right)}{k}} }{1- e^{-2 \pi n v}} \nonumber \\ 
 &-\frac{\pi i}{2} +3 \pi i s(h,k) -\frac{\pi }{4k} (v- \frac{1}{v}) + \frac{\pi z^2 k}{v} + \pi i z - \frac{\pi z}{v} = -\frac{1}{2} \log(v).
\end{align}
Now we introduce the following function and use the residue theorem on contour C shown in figure 1 to prove (14) under the condition $v$ is real and positive and $v> |\Im(z)| > 0$. Then we fix $z$, such that $\Im(z) \ne 0$ and $z$ is not a zero of $\vartheta_1(z, \tau) $ i.e $z\neq m+n\tau$, and use analytic continuation to extend $v$ to the full plane. Once done, we extend $z$ to the full plane via analytic continuation. \\ \\ For purpose of symmetry, we introduce the variable $w$, such that $w\equiv Hh~(mod~k)$ and $1 \leq w\leq k-1$. \\ \\
Now for $\displaystyle m \in \mathbb{N} $ and $ \displaystyle~ N= m +\frac{1}{2}$, define
\begin{align}
 F_N(x) = & \frac{1}{4ix} \mathrm{\coth (\pi N \textit{x}) \, \cot(\pi N \textit{xv})} + \sum_{\mu=1}^{k-1} A_N(x) + 2\sum_{\mu=1}^{k-1} e^{2 \pi N z x} A_N(x) \nonumber \\ &+ \frac{e^{2 \pi Nz x}}{x} \left(\frac{1}{1-e^{2 \pi N x}} \right) \left( \frac{e^{2 \pi i Nv x}}{1- e^{2 \pi i Nv x}} \right) + \frac{e^{-2 \pi Nz x}}{x} \left(\frac{e^{2 \pi N x}}{1-e^{2 \pi N x}} \right) \left( \frac{1}{1- e^{2 \pi i Nv x}} \right) , \end{align}
where 
\begin{align}
 A_N(x) = \frac{1}{x} \cdot \frac{e^{2 \pi N w  x/k}}{1-e^{2\pi N x}} \cdot \frac{e^{2 \pi i N  \mu v x/k}}{1-e^{-2 \pi i   N x v}}. \end{align}
This function has a pole of order 3 at 0 and two simple poles at $\displaystyle \frac{in}{N}$ and $\displaystyle -\frac{n}{Nv} .$ \\ \\
We will find the residues at each pole for each separate sub-function of  $F_N(x)$ and then group them to get the residue of $F_N(x)$.\\ \\
\underline{At $\displaystyle x= 0$} \\ \\
A straightforward calculation shows that the residue at 0 of \[\frac{1}{4ix}\mathrm{\coth(\pi \textit{Nx}) 
\, \cot(\pi \textit{Nvx})}\] is \[\frac{i}{12}\left(v-\frac{1}{v}\right).\]
 And the residue at zero of  \[A_(x)= \frac{1}{x} \cdot \frac{e^{2 \pi N w  x/k}}{1-e^{2\pi N x}} \cdot \frac{e^{2 \pi i N  \mu v x/k}}{1-e^{-2 \pi i   N x v}}\] is 
\begin{align} \left(\frac{1}{12}-\frac{\mu}{2k}+\frac{1}{2}\frac{\mu^2}{k^2}\right)vi +\left(\frac{\mu}{k}-\frac{1}{2}\right)\left(\frac{w}{k}-\frac{1}{2}\right)+\left(\frac{1}{12}-\frac{w}{2k}+\frac{w^2}{2k^2}\right)\frac{1}{vi}.
\end{align}
To get the residue of the sub-function \[\sum_{\mu=1}^{k-1} A_n(x),\] we sum over $\mu$ from 1 to $k-1$, however since $ w\equiv hu~ \pmod{k}$ then $w$ also runs from 1 to $k-1$.\\ \\
We also note that by definition as in \cite{Rd}
\[\sum_{\mu=1}^{k-1}\left(\frac{\mu}{k}-\frac{1}{2}\right)\left(\frac{w}{k}-\frac{1}{2}\right)=s(h,k).\]
Hence summing (17) from 1 to $k-1$ we have
\begin{align} 
\left(-\frac{1}{12}+\frac{1}{12k}\right)vi +s(h,k) +\left(-\frac{1}{12}+\frac{1}{12k}\right)\frac{1}{vi}.
\end{align}
Now for 
\begin{align}
2\sum_{\mu=1}^{k-1} e^{2 \pi N z x} A_n(x),
\end{align}
we use the Taylor expansion of $e$ to calculate the residue at 0
\[e^{2\pi Nzx}=1+2\pi Nz+...\]
We get that (12) yields twice (11) since the sum is multiplied by a factor of 2
\begin{align} 
2\left(-\frac{1}{12}+\frac{1}{12k}\right)vi +2s(h,k) +2\left(-\frac{1}{12}+\frac{1}{12k}\right)\frac{1}{vi}
\end{align}
in addition to the following terms containing $z$ summed over 1 to $k-1$
\[\sum_{\mu=1}^{k-1}\left(-z-\frac{z}{iv}+\frac{z^2}{iv}\right)+\sum_{\mu=1}^{k-1} \frac{2zw}{ivk}+\sum_{\mu=1}^{k-1}\frac{2z\mu}{k}, \]
which yields
\begin{align}
\frac{z^2}{iv}(k-1).
\end{align}
Finally, the residue at 0 of
\[\frac{e^{2 \pi Nz x}}{x} \left(\frac{1}{1-e^{2 \pi N x}} \right) \left( \frac{e^{2 \pi i Nv x}}{1- e^{2 \pi i Nv x}} \right) + \frac{e^{-2 \pi Nz x}}{x} \left(\frac{e^{2 \pi N x}}{1-e^{2 \pi N x}} \right) \left( \frac{1}{1- e^{2 \pi i Nv x}} \right)\]
is
\[-\frac{z+iv}{2} -\frac{z+iv}{2vi} +\frac{z^2+2zvi-v^2}{2vi}+\frac{1}{4}+\frac{z-1}{2}+\frac{z-1}{2vi}+\frac{z^2-2z+1}{2vi}+\frac{1}{4}\]
which eventually yields
\begin{align}
    \frac{z^2}{vi}-\frac{z}{vi}+z-\frac{1}{2}.
\end{align}
Summing (18), (20), (21) and (22), i.e the residue of the whole fuction $F_N(x)$ at 0, we obtain the final result
\begin{align}
    \frac{z^2k}{vi}-\frac{z}{vi}+z+\frac{i}{4k}\left(v-\frac{1}{v}\right) +3s(h,k).
\end{align}
For the residues of $F_N(x)$ at the simple poles, we proceed similarly. \\ \\ \\ \\
\underline{At  $\displaystyle x= \frac{in}{N}$ } \\ \\
 The residue of 
 \begin{align*}
  \frac{1}{4ix} \cdot \coth (\pi N x) \cdot \cot(\pi N x v)  \end{align*} is \begin{align*} \frac{i}{4 \pi  n} coth( \pi n v) =\frac{i}{4 \pi n }   \left(1+ \frac{2 e^{-2 \pi n v}}{1- e^{-2 \pi n v} } \right) .\end{align*}
 By the definition of the variable $w$, \\ \\
 $ \displaystyle
 Res \left(  \sum_{\mu=1}^{k-1} A_N(x), \frac{in}{N} \right) =  - \frac{1}{2 \pi i} \sum_{\mu=1}^{k-1} \frac{1}{n} e^{\frac{2\pi i n w}{k}} \frac{e^{\frac{-2\pi n v  \mu}{k}}}{1- e^{-2 \pi n v}}=  - \frac{1}{2 \pi i} \sum_{\mu=1}^{k-1} \frac{1}{n} e^{\frac{2\pi i n h  \mu}{k}} \frac{e^{\frac{-2\pi n v  \mu}{k}}}{1- e^{-2 \pi n v}}.\\ \\$
 Therefore, 
  \[ \displaystyle \mathrm{Res} \left( \sum_{\mu=1}^{k-1} e^{2 \pi N z x} A_N(x), \frac{in}{N} \right) =   - \frac{1}{2 \pi i}\sum_{\mu=1}^{k-1}\frac{1}{n} e^{2 \pi i n z} e^{\frac{2\pi i n h  \mu}{k-1}} \frac{e^{\frac{-2\pi n v  \mu}{k}}}{1- e^{-2 \pi n v}} .\]
The residue of the last two terms in the function is   
\[ \displaystyle 
-\frac{1}{2 \pi i n} e^{2 \pi i n z} \frac{e^{-2 \pi nv}}{1- e^{-2 \pi n v}} -\frac{1}{2 \pi i n} e^{-2 \pi i n z} \frac{1}{1- e^{-2 \pi n v}} .\]
The parallelogram C contains the poles   $\displaystyle x= \frac{in}{N}$  for $ \displaystyle -m \leq n \leq -1 $ and $\displaystyle 1 \leq n \leq m, $ so combining the sub-results, summing over these poles, and re-indexing the sum over the negative n, we obtain \\ \\
$ \displaystyle  \frac{i}{2 \pi  } \sum_{n=1}^{m} \frac{1}{n} \left(1+ \frac{2 e^{-2 \pi n v}}{1- e^{-2 \pi n v} } \right) - \frac{1}{ \pi i } \sum_{n=1}^{m} \sum_{\mu=1}^{k-1} \frac{1}{n} e^{\frac{2\pi i n h  \mu}{k}} \frac{e^{\frac{-2\pi n v  \mu}{k}}}{1- e^{-2 \pi n v}}  \\ \\ - \frac{1}{ \pi i } \sum_{n=1}^{m} \sum_{\mu=1}^{k-1}\frac{1}{n} e^{2 \pi i n z} e^{\frac{2\pi i n h  \mu}{k}} \frac{e^{\frac{-2\pi n v  \mu}{k}}}{1- e^{-2 \pi n v}} 
- \frac{1}{ \pi i } \sum_{n=1}^{m} \sum_{\mu=1}^{k-1}\frac{1}{n} e^{-2 \pi i n z} e^{\frac{2\pi i n h  \mu}{k}} \frac{e^{\frac{-2\pi n v  \mu}{k}}}{1- e^{-2 \pi n v}} \\ \\
- \frac{1}{ \pi i } \sum_{n=1}^{m} \frac{1}{n} e^{2 \pi i n z} \frac{e^{-2\pi n v  }}{1- e^{-2 \pi n v}} - \frac{1}{ \pi i } \sum_{n=1}^{m} \frac{1}{n} e^{-2 \pi i n z} \frac{1}{1- e^{-2 \pi n v}}.
$
\\  \\ \\ \\ \\ \\
The first summation includes the $k^{th}$ term of the first double sum, and the last two summations match the $k^{th}$ and $0^{th}$ terms of the second and third double sums respectively, so we get \\ \\
 $ \displaystyle  \frac{i}{2 \pi  } \sum_{n=1}^{m} \frac{1}{n} - \frac{1}{ \pi i } \sum_{n=1}^{m} \sum_{\mu=1}^{k} \frac{1}{n} e^{\frac{2\pi i n h  \mu}{k}} \frac{e^{\frac{-2\pi n v  \mu}{k}}}{1- e^{-2 \pi n v}} - \frac{1}{ \pi i } \sum_{n=1}^{m} \sum_{\mu=1}^{k}\frac{1}{n} e^{2 \pi i n z} e^{\frac{2\pi i n h  \mu}{k}} \frac{e^{\frac{-2\pi n v  \mu}{k}}}{1- e^{-2 \pi n v}} \\ \\
- \frac{1}{ \pi i } \sum_{n=1}^{m} \sum_{\mu=0}^{k-1}\frac{1}{n} e^{-2 \pi i n z} e^{\frac{2\pi i n h  \mu}{k}} \frac{e^{\frac{-2\pi n v  \mu}{k}}}{1- e^{-2 \pi n v}}. ~~~~~~~~~~~~~~ (24)
$\\ 
\vspace{0.43cm}
\\
 \underline{At $x=\displaystyle \frac{-n}{Nv}$}\\ \\
 The residue at $x=\displaystyle \frac{-n}{Nv}$  of \[\displaystyle \frac{1}{4ix} \cdot \coth (\pi N x) \cdot \cot(\pi N x v) \] is  \[ \displaystyle \frac{i}{4 \pi  n} coth( - \frac{\pi n}{v})= \frac{1}{4 \pi i n }  (1+ \frac{2 e^{ -\frac{2 \pi n}{v}}}{1- e^{-\frac{2 \pi n }{v}} }). \] 
 And,
 \[\displaystyle Res \left(  \sum_{\mu=1}^{k-1} A_N(x), \frac{in}{N} \right) =  \frac{1}{2 \pi i}\sum_{\mu=1}^{k-1} \frac{1}{n} e^{\frac{-2\pi  n w}{kv}} \frac{e^{\frac{-2\pi i n   \mu}{k}}}{1- e^{-2 \pi n/ v}}.\]
 But by our choice of $Hh\equiv -1~ \pmod{k}$, we have
  \[ \displaystyle Res \left(  \sum_{\mu=1}^{k-1} A_N(x), \frac{in}{N} \right) =  \frac{1}{2 \pi i}\sum_{\mu=1}^{k-1} \frac{1}{n} e^{\frac{-2\pi  n w}{kv}} \frac{e^{\frac{2\pi i n H h  \mu}{k}}}{1- e^{-2 \pi n/ v}}.\]
 Thus, by our choice of the variable $w$, we have
\[\displaystyle Res \left(  \sum_{\mu=1}^{k-1} A_N(x), \frac{in}{N} \right) = \frac{1}{2 \pi i}\sum_{w=1}^{k-1} \frac{1}{n} e^{\frac{-2\pi  n w}{kv}} \frac{e^{\frac{2\pi i n  H w}{k}}}{1- e^{-2 \pi n/ v}}.\]
Similarly, we get
\[\displaystyle Res \left( \sum_{\mu=1}^{k-1} e^{2 \pi N z x} A_N(x), \frac{in}{N} \right) = \sum_{w=1}^{k-1} \frac{1}{n} e^{2 \pi  n \frac{z}{v}} e^{\frac{-2\pi  n w}{kv}} \frac{e^{\frac{2\pi i n  H w}{k}}}{1- e^{-2 \pi n/ v}}. \]
The parallelogram C contains the poles   $\displaystyle x= \frac{-n}{Nv}$  for $ \displaystyle -m \leq n \leq -1 $ and $\displaystyle 1 \leq n \leq m $ , so combining the sub-results, summing over these poles and re-indexing the sum over the negative n, we obtain \\ \\ \\
$ \displaystyle  \frac{1}{2 \pi i } \sum_{n=1}^{m} \frac{1}{n} (1+ \frac{2 e^{- 2 \pi n /v} }{1- e^{-2 \pi n /v} }) + \frac{1}{ \pi i } \sum_{n=1}^{m} \sum_{\mu=1}^{k-1} \frac{1}{n} e^{\frac{-2\pi  n w}{kv}} \frac{e^{\frac{2\pi i n  H w}{k}}}{1- e^{-2 \pi n/ v}} \\ \\+ \frac{1}{ \pi i } \sum_{n=1}^{m} \sum_{\mu=1}^{k-1}\frac{1}{n} e^{2 \pi i n z}  e^{\frac{-2\pi  n w}{kv}} \frac{e^{\frac{2\pi i n  H w}{k}}}{1- e^{-2 \pi n/ v}} 
+ \frac{1}{ \pi i } \sum_{n=1}^{m} \sum_{\mu=1}^{k-1}\frac{1}{n} e^{-2 \pi i n z}  e^{\frac{-2\pi  n w}{kv}} \frac{e^{\frac{2\pi i n  H w}{k}}}{1- e^{-2 \pi n/ v}} \\ \\ 
+ \frac{1}{ \pi i } \sum_{n=1}^{m} \frac{1}{n} e^{-2 \pi  n z/v} \frac{1}{1- e^{-2 \pi n/ v}} + \frac{1}{ \pi i } \sum_{n=1}^{m} \frac{1}{n} e^{2 \pi  n z/v} \frac{e^{-2\pi n/  v  }}{1- e^{-2 \pi n/ v}}.
$
\\ \\
The first summation includes the $k^{th}$ term of the first double sum, and the last two summations match the $0^{th}$ and $k^{th}$ terms of the third and second double sums respectively, so we get \\ \\ \\
$ \displaystyle  \frac{1}{2 \pi i } \sum_{n=1}^{m} \frac{1}{n} + \frac{1}{ \pi i } \sum_{n=1}^{m} \sum_{\mu=1}^{k} \frac{1}{n} e^{\frac{-2\pi  n w}{kv}} \frac{e^{\frac{2\pi i n  H w}{k}}}{1- e^{-2 \pi n/ v}}+ \frac{1}{ \pi i } \sum_{n=1}^{m} \sum_{\mu=1}^{k}\frac{1}{n} e^{2 \pi i n z}  e^{\frac{-2\pi  n w}{kv}} \frac{e^{\frac{2\pi i n  H w}{k}}}{1- e^{-2 \pi n/ v}} \\ \\
+ \frac{1}{ \pi i } \sum_{n=1}^{m} \sum_{\mu=0}^{k-1}\frac{1}{n} e^{-2 \pi i n z} \frac{1}{n} e^{\frac{-2\pi  n w}{kv}} \frac{e^{\frac{2\pi i n  H w}{k}}}{1- e^{-2 \pi n/ v}}. ~~~~~~~~~~~~~~ (25)
$
\\ \\
Now note that if we multiply the sum of all of the residues of $F_N(x)$ by $\pi i$, i.e. $ \displaystyle \pi i~ \{ (23) + (24) + (25)  \}$ and take the limit m $\rightarrow \infty$, we precisely get the left-hand side of (14).  \\ \\
What remains to show to finish the proof of (14) is that 
\begin{align*}
\displaystyle
\frac{1}{2}\lim_{m \to \infty} \int_{C} F_N(x) d \zeta = - \frac{1}{2}\log(v).
\end{align*}  
One can easily show that \(
\lim_{m \to \infty} xF_N(x)= -\frac{1}{4}\)  on edges $(\frac{1}{v},i), (-\frac{1}{v}, -i)$ and   $ \frac{1}{4}$ on the other two edges, and  $F_N(x)$ is uniformly bounded on C. Thus by appealing to Arzela's bounded convergence theorem, we obtain
\\
\begin{align*}
\lim_{m \to \infty} \int_{C} F_N(x) d x &= \lim_{m \to \infty} \int_{C} xF_N(x) \frac{d x}{x}  = \frac{1}{4} \left( -\int_{\frac{1}{v}}^{i} + \int_{i}^{-\frac{1}{v}} -\int_{-\frac{1}{v}}^{-i} +\int_{-i}^{\frac{1}{v}} \right)\frac{d x}{x}\\
&= \frac{1}{2} \left( - \int_{\frac{1}{v}}^{i} +\int_{-i}^{\frac{1}{v}} \right)\frac{d x}{x}= \frac{1}{2} \left(-\frac{\pi i}{2} +\log(\frac{1}{v}) +\log(\frac{1}{v})  + \frac{\pi i}{2} \right) \\
&= -\log(v). 
\end{align*}
This completes the proof of (14) and therefore  the theorem.\\
\begin{center}
ACKNOWLEDGEMENT
\end{center} 
We would like to express our deepest gratitude to our advisor Professor Wissam Raji along with the Department of Mathematics at the American University of Beirut (AUB). We also like to thank the Center of Advanced Mathematical Sciences
(CAMS) at AUB for the guidance and support we received at the summer research camp (SRC). Additionally, we would like to thank the referee whose valuable comments helped us enhance the content and the structure of our paper.\\

\vspace{0.25cm}
Department of Mathematics, American University of Beirut, Beirut, Lebanon
\textit{E-mail address: mmm133@mail.aub.edu} \\
\textit{E-mail address: ays11@mail.aub.edu}

\begin{thebibliography}{unsrt}
\bibitem{Ap} Apostol, T. Modular Functions and Dirichlet Series in Number Theory,
Springer-Verlag, New York, (1989). MR1027834 (52-190)
\bibitem{dobholkar2012quantum}
A. Dabholkar, S. Murthy, and D. Zagier
\bibitem{g. iseki}
M. Me'meh \& A. Saraeb. A Generalization of Iseki’s Formula and The Jacobi
Theta Function. Hardy-Ramanujan Journal. vol. 45, no 10, (2022), pp. 130-139. 
\bibitem{g. gordon}
M. Me'meh \& A. Saraeb. A New Proof for the Transformation Laws of The Jacobi Theta functions. Submitted for publication.\\
https://arxiv.org/abs/2207.12806
\bibitem{seigelGen}
Me'meh, Maher and Saraeb, Ali. On a New Generalization of Siegel's Method to $\theta_1(z,\tau)$ , submitted for publication, https://arxiv.org/abs/2208.02168.
\bibitem{Rd}
H.  Rademacher, On the transformation of logn(r). J. Indian Math Soc. 19(1955), 25-30. MR0070660 (17:15f)
\bibitem{Raji}
  Raji, Wissam. A New Proof of the Transformation Law of Jacobi's Theta Function $\vartheta_3(w,\tau)$. Proceedings of the American Mathematical Society, vol. 135, no. 10, (2007), pp. 3127–32,
 http://www.jstor.org/stable/20534932,(3127-3132)
\bibitem{Rich}
Richter, Olav K. “On The Transformation Laws For Theta Functions.” The Rocky Mountain Journal of Mathematics, vol. 34, no. 4, 2004, pp. 1473–81. JSTOR, http://www.jstor.org/stable/44239039. 
\bibitem{Yurii}
Brezhnev, Yurii V. “Non-Canonical Extension of $\vartheta$-Functions and Modular Integrability of $\vartheta$-Constants.” Proceedings of the Royal Society of Edinburgh: Section A Mathematics, vol. 143, no. 4, 2013, pp. 689–738., doi:10.1017/S0308210512001023. (709-710)

\end{thebibliography}
\end{document}